\pdfoutput=1

\documentclass[12pt,a4paper]{article}

\usepackage{amsfonts,amssymb,amsthm,amsmath,latexsym}
\usepackage{indent}
\usepackage{subeqnarray, eepicemu, url,cite,bm}
\usepackage{multirow}  
\usepackage{tensor}  
\usepackage{graphicx,graphbox}
\usepackage{float}  

\date{April 3, 2022}

\begin{document}

\title{\vspace*{-1cm}
       When does a hypergeometric function ${}_{p\!}F_q$ \\[2mm]
       belong to the Laguerre--P\'olya class $LP^+$?
      }

\author{
     {\small Alan D.~Sokal}                  \\[2mm]
     {\small\it Department of Mathematics}   \\[-2mm]
     {\small\it University College London}   \\[-2mm]
     {\small\it Gower Street}                \\[-2mm]
     {\small\it London WC1E 6BT}             \\[-2mm]
     {\small\it UNITED KINGDOM}              \\[-2mm]
     {\small\tt sokal@math.ucl.ac.uk}        \\[-2mm]
     {\protect\makebox[5in]{\quad}}  
     \\[-2mm]
     {\small\it Department of Physics}       \\[-2mm]
     {\small\it New York University}         \\[-2mm]
     {\small\it 726 Broadway}                \\[-2mm]
     {\small\it New York, NY 10003}          \\[-2mm]
     {\small\it USA}                         \\[-2mm]
     {\small\tt sokal@nyu.edu}               \\[3mm]
}

\maketitle
\thispagestyle{empty}   

\begin{abstract}
I show that a hypergeometric function
$\tensor[_{p \!}]{F}{_{q}}\!(a_1,\ldots,a_p;b_1,\ldots,b_q;\,\cdot\,)$
with ${p \le q}$ belongs to the Laguerre--P\'olya class $LP^+$
for arbitrarily large $b_{p+1},\ldots,b_q > 0$
if and only if, after a possible reordering,
the differences $a_i - b_i$ are nonnegative integers.
This result arises as an easy corollary of the case $p=q$
proven two decades ago by Ki and Kim.
I also give explicit examples for the case
$\tensor[_{1 \!}]{F}{_{2}}\!$.
\end{abstract}


\bigskip
\bigskip
\noindent
{\bf Mathematics Subject Classification (MSC 2010) codes:}
33C20 (Primary);
30B70, 30C15, 30D20, 30E05, 44A60 (Secondary).

\clearpage


\newtheorem{theorem}{Theorem}
\newtheorem{proposition}[theorem]{Proposition}
\newtheorem{lemma}[theorem]{Lemma}
\newtheorem{corollary}[theorem]{Corollary}
\newtheorem{definition}[theorem]{Definition}
\newtheorem{conjecture}[theorem]{Conjecture}
\newtheorem{question}[theorem]{Question}
\newtheorem{problem}[theorem]{Problem}
\newtheorem{openproblem}[theorem]{Open Problem}
\newtheorem{example}[theorem]{Example}

\renewcommand{\theenumi}{\alph{enumi}}
\renewcommand{\labelenumi}{(\theenumi)}
\def\eop{\hbox{\kern1pt\vrule height6pt width4pt
depth1pt\kern1pt}\medskip}
\def\prf{\par\noindent{\bf Proof.\enspace}\rm}
\def\rmk{\par\medskip\noindent{\bf Remark\enspace}\rm}

\newcommand{\textbfit}[1]{\textbf{\textit{#1}}}

\newcommand{\bigdash}{%
\smallskip\begin{center} \rule{5cm}{0.1mm} \end{center}\smallskip}

\newcommand{\safepar}{ {\protect\hfill\protect\break\hspace*{5mm}} }

\newcommand{\be}{\begin{equation}}
\newcommand{\ee}{\end{equation}}
\newcommand{\<}{\langle}
\renewcommand{\>}{\rangle}
\newcommand{\widebar}{\overline}
\def\reff#1{(\protect\ref{#1})}
\def\spose#1{\hbox to 0pt{#1\hss}}
\def\ltapprox{\mathrel{\spose{\lower 3pt\hbox{$\mathchar"218$}}
    \raise 2.0pt\hbox{$\mathchar"13C$}}}
\def\gtapprox{\mathrel{\spose{\lower 3pt\hbox{$\mathchar"218$}}
    \raise 2.0pt\hbox{$\mathchar"13E$}}}
\def\textprime{${}^\prime$}
\def\proof{\par\medskip\noindent{\sc Proof.\ }}
\def\firstproof{\par\medskip\noindent{\sc First Proof.\ }}
\def\secondproof{\par\medskip\noindent{\sc Second Proof.\ }}
\def\alternateproof{\par\medskip\noindent{\sc Alternate Proof.\ }}
\def\algebraicproof{\par\medskip\noindent{\sc Algebraic Proof.\ }}
\def\combinatorialproof{\par\medskip\noindent{\sc Combinatorial Proof.\ }}
\def\proofof#1{\bigskip\noindent{\sc Proof of #1.\ }}
\def\firstproofof#1{\bigskip\noindent{\sc First Proof of #1.\ }}
\def\secondproofof#1{\bigskip\noindent{\sc Second Proof of #1.\ }}
\def\thirdproofof#1{\bigskip\noindent{\sc Third Proof of #1.\ }}
\def\algebraicproofof#1{\bigskip\noindent{\sc Algebraic Proof of #1.\ }}
\def\combinatorialproofof#1{\bigskip\noindent{\sc Combinatorial Proof of #1.\ }}
\def\sketchofproof{\par\medskip\noindent{\sc Sketch of proof.\ }}
\renewcommand{\qed}{ $\square$ \bigskip}
\newcommand{\myendremark}{ $\blacksquare$ \bigskip}
\def\half{ {1 \over 2} }
\def\third{ {1 \over 3} }
\def\twothird{ {2 \over 3} }
\def\smfrac#1#2{{\textstyle{#1\over #2}}}
\def\smhalf{ {\smfrac{1}{2}} }
\newcommand{\real}{\mathop{\rm Re}\nolimits}
\renewcommand{\Re}{\mathop{\rm Re}\nolimits}
\newcommand{\imag}{\mathop{\rm Im}\nolimits}
\renewcommand{\Im}{\mathop{\rm Im}\nolimits}
\newcommand{\sgn}{\mathop{\rm sgn}\nolimits}
\newcommand{\tr}{\mathop{\rm tr}\nolimits}
\newcommand{\supp}{\mathop{\rm supp}\nolimits}
\newcommand{\disc}{\mathop{\rm disc}\nolimits}
\newcommand{\diag}{\mathop{\rm diag}\nolimits}
\newcommand{\tridiag}{\mathop{\rm tridiag}\nolimits}
\newcommand{\AZ}{\mathop{\rm AZ}\nolimits}
\newcommand{\NC}{\mathop{\rm NC}\nolimits}
\newcommand{\PF}{{\rm PF}}
\newcommand{\rk}{\mathop{\rm rk}\nolimits}
\newcommand{\perm}{\mathop{\rm perm}\nolimits}
\def\hboxscript#1{ {\hbox{\scriptsize\em #1}} }
\renewcommand{\emptyset}{\varnothing}
\newcommand{\eqdef}{\stackrel{\rm def}{=}}

\newcommand{\restrict}{\upharpoonright}

\newcommand{\compinv}{{\langle -1 \rangle}}   

\newcommand{\divides}{\mid}
\newcommand{\notdivides}{\nmid}

\newcommand{\scra}{{\mathcal{A}}}
\newcommand{\scrb}{{\mathcal{B}}}
\newcommand{\scrc}{{\mathcal{C}}}
\newcommand{\scrd}{{\mathcal{D}}}
\newcommand{\scrdtilde}{{\widetilde{\mathcal{D}}}}
\newcommand{\scre}{{\mathcal{E}}}
\newcommand{\scrf}{{\mathcal{F}}}
\newcommand{\scrg}{{\mathcal{G}}}
\newcommand{\scrh}{{\mathcal{H}}}
\newcommand{\scri}{{\mathcal{I}}}
\newcommand{\scrj}{{\mathcal{J}}}
\newcommand{\scrk}{{\mathcal{K}}}
\newcommand{\scrl}{{\mathcal{L}}}
\newcommand{\scrm}{{\mathcal{M}}}
\newcommand{\scrn}{{\mathcal{N}}}
\newcommand{\scro}{{\mathcal{O}}}
\newcommand\scroo{
  \mathchoice
    {{\scriptstyle\mathcal{O}}}
    {{\scriptstyle\mathcal{O}}}
    {{\scriptscriptstyle\mathcal{O}}}
    {\scalebox{0.6}{$\scriptscriptstyle\mathcal{O}$}}
  }
\newcommand{\scrp}{{\mathcal{P}}}
\newcommand{\scrq}{{\mathcal{Q}}}
\newcommand{\scrr}{{\mathcal{R}}}
\newcommand{\scrs}{{\mathcal{S}}}
\newcommand{\scrt}{{\mathcal{T}}}
\newcommand{\scrv}{{\mathcal{V}}}
\newcommand{\scrw}{{\mathcal{W}}}
\newcommand{\scrz}{{\mathcal{Z}}}
\newcommand{\SP}{{\mathcal{SP}}}
\newcommand{\ST}{{\mathcal{ST}}}

\newcommand{\bfa}{{\mathbf{a}}}
\newcommand{\bfb}{{\mathbf{b}}}
\newcommand{\bfc}{{\mathbf{c}}}
\newcommand{\bfd}{{\mathbf{d}}}
\newcommand{\bfe}{{\mathbf{e}}}
\newcommand{\bfh}{{\mathbf{h}}}
\newcommand{\bfj}{{\mathbf{j}}}
\newcommand{\bfi}{{\mathbf{i}}}
\newcommand{\bfk}{{\mathbf{k}}}
\newcommand{\bfl}{{\mathbf{l}}}
\newcommand{\bfL}{{\mathbf{L}}}
\newcommand{\bfm}{{\mathbf{m}}}
\newcommand{\bfn}{{\mathbf{n}}}
\newcommand{\bfp}{{\mathbf{p}}}
\newcommand{\bfP}{{\mathbf{P}}}
\newcommand{\bfr}{{\mathbf{r}}}
\newcommand{\bfs}{{\mathbf{s}}}
\newcommand{\bfu}{{\mathbf{u}}}
\newcommand{\bfv}{{\mathbf{v}}}
\newcommand{\bfw}{{\mathbf{w}}}
\newcommand{\bfx}{{\mathbf{x}}}
\newcommand{\bfX}{{\mathbf{X}}}
\newcommand{\bfy}{{\mathbf{y}}}
\newcommand{\bfz}{{\mathbf{z}}}
\renewcommand{\k}{{\mathbf{k}}}
\newcommand{\n}{{\mathbf{n}}}
\newcommand{\vv}{{\mathbf{v}}}
\newcommand{\bv}{{\mathbf{v}}}
\newcommand{\w}{{\mathbf{w}}}
\newcommand{\x}{{\mathbf{x}}}
\newcommand{\y}{{\mathbf{y}}}
\newcommand{\cc}{{\mathbf{c}}}
\newcommand{\zero}{{\mathbf{0}}}
\newcommand{\one}{{\mathbf{1}}}
\newcommand{\bmm}{{\mathbf{m}}}

\newcommand{\ahat}{{\widehat{a}}}
\newcommand{\Zhat}{{\widehat{Z}}}

\newcommand{\C}{{\mathbb C}}
\newcommand{\D}{{\mathbb D}}
\newcommand{\Z}{{\mathbb Z}}
\newcommand{\N}{{\mathbb N}}
\newcommand{\Q}{{\mathbb Q}}
\newcommand{\PP}{{\mathbb P}}
\newcommand{\R}{{\mathbb R}}
\newcommand{\RR}{{\mathbb R}}
\newcommand{\E}{{\mathbb E}}

\newcommand{\Sym}{{\mathfrak{S}}}
\newcommand{\SymB}{{\mathfrak{B}}}
\newcommand{\Alt}{{\mathrm{Alt}}}

\newcommand{\germanA}{{\mathfrak{A}}}
\newcommand{\germanB}{{\mathfrak{B}}}
\newcommand{\germanQ}{{\mathfrak{Q}}}
\newcommand{\germanh}{{\mathfrak{h}}}

\newcommand{\myle}{\preceq}
\newcommand{\myge}{\succeq}
\newcommand{\mygt}{\succ}

\newcommand{\B}{{\sf B}}
\newcommand{\OB}{B^{\rm ord}}
\newcommand{\OS}{{\sf OS}}
\newcommand{\OO}{{\sf O}}
\newcommand{\OSP}{{\sf OSP}}
\newcommand{\Eu}{{\sf Eu}}
\newcommand{\ERR}{{\sf ERR}}
\newcommand{\sfB}{{\sf B}}
\newcommand{\sfD}{{\sf D}}
\newcommand{\sfE}{{\sf E}}
\newcommand{\sfG}{{\sf G}}
\newcommand{\sfJ}{{\sf J}}
\newcommand{\sfL}{{\sf L}}
\newcommand{\sfLhat}{{\widehat{{\sf L}}}}
\newcommand{\sfLtilde}{{\widetilde{{\sf L}}}}
\newcommand{\sfP}{{\sf P}}
\newcommand{\sfQ}{{\sf Q}}
\newcommand{\sfS}{{\sf S}}
\newcommand{\sfT}{{\sf T}}
\newcommand{\sfW}{{\sf W}}
\newcommand{\sfMV}{{\sf MV}}
\newcommand{\AMV}{{\sf AMV}}
\newcommand{\BM}{{\sf BM}}
\newcommand{\emIB}{B^{\rm irr}}
\newcommand{\emIP}{P^{\rm irr}}
\newcommand{\emOB}{B^{\rm ord}}
\newcommand{\emCB}{B^{\rm cyc}}
\newcommand{\emSC}{P^{\rm cyc}}

\newcommand{\lev}{{\rm lev}}
\newcommand{\stat}{{\rm stat}}
\newcommand{\cyc}{{\rm cyc}}
\newcommand{\mysteryone}{{\rm mys1}}
\newcommand{\mysterytwo}{{\rm mys2}}
\newcommand{\Asc}{{\rm Asc}}
\newcommand{\asc}{{\rm asc}}
\newcommand{\Des}{{\rm Des}}
\newcommand{\des}{{\rm des}}
\newcommand{\Exc}{{\rm Exc}}
\newcommand{\exc}{{\rm exc}}
\newcommand{\Wex}{{\rm Wex}}
\newcommand{\wex}{{\rm wex}}
\newcommand{\Fix}{{\rm Fix}}
\newcommand{\fix}{{\rm fix}}
\newcommand{\lrmax}{{\rm lrmax}}
\newcommand{\rlmax}{{\rm rlmax}}
\newcommand{\Rec}{{\rm Rec}}
\newcommand{\rec}{{\rm rec}}
\newcommand{\Arec}{{\rm Arec}}
\newcommand{\arec}{{\rm arec}}
\newcommand{\ERec}{{\rm ERec}}
\newcommand{\erec}{{\rm erec}}
\newcommand{\EArec}{{\rm EArec}}
\newcommand{\earec}{{\rm earec}}
\newcommand{\recarec}{{\rm recarec}}
\newcommand{\nonrec}{{\rm nonrec}}
\newcommand{\Cpeak}{{\rm Cpeak}}
\newcommand{\cpeak}{{\rm cpeak}}
\newcommand{\Cval}{{\rm Cval}}
\newcommand{\cval}{{\rm cval}}
\newcommand{\Cdasc}{{\rm Cdasc}}
\newcommand{\cdasc}{{\rm cdasc}}
\newcommand{\Cddes}{{\rm Cddes}}
\newcommand{\cddes}{{\rm cddes}}
\newcommand{\cdrise}{{\rm cdrise}}
\newcommand{\cdfall}{{\rm cdfall}}
\newcommand{\Peak}{{\rm Peak}}
\newcommand{\peak}{{\rm peak}}
\newcommand{\Val}{{\rm Val}}
\newcommand{\val}{{\rm val}}
\newcommand{\Dasc}{{\rm Dasc}}
\newcommand{\dasc}{{\rm dasc}}
\newcommand{\Ddes}{{\rm Ddes}}
\newcommand{\ddes}{{\rm ddes}}
\newcommand{\inv}{{\rm inv}}
\newcommand{\maj}{{\rm maj}}
\newcommand{\rs}{{\rm rs}}
\newcommand{\cross}{{\rm cr}}
\newcommand{\crosshat}{{\widehat{\rm cr}}}
\newcommand{\nest}{{\rm ne}}
\newcommand{\rodd}{{\rm rodd}}
\newcommand{\reven}{{\rm reven}}
\newcommand{\lodd}{{\rm lodd}}
\newcommand{\leven}{{\rm leven}}
\newcommand{\sg}{{\rm sg}}
\newcommand{\bl}{{\rm bl}}
\newcommand{\tran}{{\rm tr}}
\newcommand{\area}{{\rm area}}
\newcommand{\ret}{{\rm ret}}
\newcommand{\peaks}{{\rm peaks}}
\newcommand{\hl}{{\rm hl}}
\newcommand{\sll}{{\rm sl}}
\newcommand{\negg}{{\rm neg}}
\newcommand{\imp}{{\rm imp}}
\newcommand{\osg}{{\rm osg}}
\newcommand{\ons}{{\rm ons}}
\newcommand{\isg}{{\rm isg}}
\newcommand{\ins}{{\rm ins}}
\newcommand{\LL}{{\rm LL}}
\newcommand{\height}{{\rm ht}}
\newcommand{\as}{{\rm as}}

\newcommand{\ba}{{\bm{a}}}
\newcommand{\bahat}{{\widehat{\bm{a}}}}
\newcommand{\sfa}{{{\sf a}}}
\newcommand{\bb}{{\bm{b}}}
\newcommand{\bc}{{\bm{c}}}
\newcommand{\bchat}{{\widehat{\bm{c}}}}
\newcommand{\bd}{{\bm{d}}}
\newcommand{\bee}{{\bm{e}}}
\newcommand{\beh}{{\bm{eh}}}
\newcommand{\bff}{{\bm{f}}}
\newcommand{\bg}{{\bm{g}}}
\newcommand{\bh}{{\bm{h}}}
\newcommand{\bll}{{\bm{\ell}}}
\newcommand{\bp}{{\bm{p}}}
\newcommand{\br}{{\bm{r}}}
\newcommand{\bs}{{\bm{s}}}
\newcommand{\bu}{{\bm{u}}}
\newcommand{\bw}{{\bm{w}}}
\newcommand{\bx}{{\bm{x}}}
\newcommand{\by}{{\bm{y}}}
\newcommand{\bz}{{\bm{z}}}
\newcommand{\bA}{{\bm{A}}}
\newcommand{\bB}{{\bm{B}}}
\newcommand{\bC}{{\bm{C}}}
\newcommand{\bE}{{\bm{E}}}
\newcommand{\bF}{{\bm{F}}}
\newcommand{\bG}{{\bm{G}}}
\newcommand{\bH}{{\bm{H}}}
\newcommand{\bI}{{\bm{I}}}
\newcommand{\bJ}{{\bm{J}}}
\newcommand{\bL}{{\bm{L}}}
\newcommand{\bLhat}{{\widehat{\bm{L}}}}
\newcommand{\bM}{{\bm{M}}}
\newcommand{\bN}{{\bm{N}}}
\newcommand{\bP}{{\bm{P}}}
\newcommand{\bQ}{{\bm{Q}}}
\newcommand{\bR}{{\bm{R}}}
\newcommand{\bS}{{\bm{S}}}
\newcommand{\bT}{{\bm{T}}}
\newcommand{\bW}{{\bm{W}}}
\newcommand{\bX}{{\bm{X}}}
\newcommand{\bY}{{\bm{Y}}}
\newcommand{\bIB}{{\bm{B}^{\rm irr}}}
\newcommand{\bOB}{{\bm{B}^{\rm ord}}}
\newcommand{\bOS}{{\bm{OS}}}
\newcommand{\bERR}{{\bm{ERR}}}
\newcommand{\bSP}{{\bm{SP}}}
\newcommand{\bMV}{{\bm{MV}}}
\newcommand{\bBM}{{\bm{BM}}}
\newcommand{\balpha}{{\bm{\alpha}}}
\newcommand{\balphapre}{{\bm{\alpha}^{\rm pre}}}
\newcommand{\bbeta}{{\bm{\beta}}}
\newcommand{\bgamma}{{\bm{\gamma}}}
\newcommand{\bGamma}{{\bm{\Gamma}}}
\newcommand{\bdelta}{{\bm{\delta}}}
\newcommand{\bkappa}{{\bm{\kappa}}}
\newcommand{\bmu}{{\bm{\mu}}}
\newcommand{\bomega}{{\bm{\omega}}}
\newcommand{\bsigma}{{\bm{\sigma}}}
\newcommand{\btau}{{\bm{\tau}}}
\newcommand{\bphi}{{\bm{\phi}}}
\newcommand{\bphihat}{{\skew{3}\widehat{\vphantom{t}\protect\smash{\bm{\phi}}}}}
\newcommand{\bpsi}{{\bm{\psi}}}
\newcommand{\bxi}{{\bm{\xi}}}
\newcommand{\bzeta}{{\bm{\zeta}}}
\newcommand{\bone}{{\bm{1}}}
\newcommand{\bzero}{{\bm{0}}}

\newcommand{\Cbar}{{\overline{C}}}
\newcommand{\Dbar}{{\overline{D}}}
\newcommand{\dbar}{{\overline{d}}}
\def\Btilde{{\widetilde{B}}}
\def\Ctilde{{\widetilde{C}}}
\def\Ftilde{{\widetilde{F}}}
\def\Gtilde{{\widetilde{G}}}
\def\Htilde{{\widetilde{H}}}
\def\Lhat{{\widehat{L}}}
\def\Ltilde{{\widetilde{L}}}
\def\Ptilde{{\widetilde{P}}}
\def\Phat{{\widehat{P}}}
\def\bfPhat{{\widehat{\bfP}}}
\def\ptilde{{\widetilde{p}}}
\def\Chat{{\widehat{C}}}
\def\ctilde{{\widetilde{c}}}
\def\zbar{{\overline{Z}}}
\def\pitilde{{\widetilde{\pi}}}
\def\omegahat{{\widehat{\omega}}}

\newcommand{\sech}{{\rm sech}}

%
%
\newcommand{\sn}{{\rm sn}}
\newcommand{\cn}{{\rm cn}}
\newcommand{\dn}{{\rm dn}}
\newcommand{\sm}{{\rm sm}}
\newcommand{\cm}{{\rm cm}}

%
%
\newcommand{\zfz}{ {{}_0 \! F_0} }
\newcommand{\zfo}{ {{}_0  F_1} }
\newcommand{\ofz}{ {{}_1 \! F_0} }
\newcommand{\ofo}{ {{}_1 \! F_1} }
\newcommand{\oft}{ {{}_1 \! F_2} }

%
%
\newcommand{\FHyper}[2]{ {\tensor[_{#1 \!}]{F}{_{#2}}\!} }
\newcommand{\FHYPER}[5]{ {\FHyper{#1}{#2} \!\biggl(
   \!\!\begin{array}{c} #3 \\[1mm] #4 \end{array}\! \bigg|\, #5 \! \biggr)} }
\newcommand{\tfo}{ {\FHyper{2}{1}} }
\newcommand{\tfz}{ {\FHyper{2}{0}} }
\newcommand{\threefz}{ {\FHyper{3}{0}} }
\newcommand{\FHYPERbottomzero}[3]{ {\FHyper{#1}{0} \hspace*{-0mm}\biggl(
   \!\!\begin{array}{c} #2 \\[1mm] \hbox{---} \end{array}\! \bigg|\, #3 \! \biggr)} }
\newcommand{\FHYPERtopzero}[3]{ {\FHyper{0}{#1} \hspace*{-0mm}\biggl(
   \!\!\begin{array}{c} \hbox{---} \\[1mm] #2 \end{array}\! \bigg|\, #3 \! \biggr)} }

\newcommand{\phiHyper}[2]{ {\tensor[_{#1}]{\phi}{_{#2}}} }
\newcommand{\psiHyper}[2]{ {\tensor[_{#1}]{\psi}{_{#2}}} }
\newcommand{\PhiHyper}[2]{ {\tensor[_{#1}]{\Phi}{_{#2}}} }
\newcommand{\PsiHyper}[2]{ {\tensor[_{#1}]{\Psi}{_{#2}}} }
\newcommand{\phiHYPER}[6]{ {\phiHyper{#1}{#2} \!\left(
   \!\!\begin{array}{c} #3 \\ #4 \end{array}\! ;\, #5, \, #6 \! \right)\!} }
\newcommand{\psiHYPER}[6]{ {\psiHyper{#1}{#2} \!\left(
   \!\!\begin{array}{c} #3 \\ #4 \end{array}\! ;\, #5, \, #6 \! \right)} }
\newcommand{\PhiHYPER}[5]{ {\PhiHyper{#1}{#2} \!\left(
   \!\!\begin{array}{c} #3 \\ #4 \end{array}\! ;\, #5 \! \right)\!} }
\newcommand{\PsiHYPER}[5]{ {\PsiHyper{#1}{#2} \!\left(
   \!\!\begin{array}{c} #3 \\ #4 \end{array}\! ;\, #5 \! \right)\!} }
\newcommand{\zerophizero}{ {\phiHyper{0}{0}} }
\newcommand{\ophizero}{ {\phiHyper{1}{0}} }
\newcommand{\zphio}{ {\phiHyper{0}{1}} }
\newcommand{\ophio}{ {\phiHyper{1}{1}} }
\newcommand{\tphio}{ {\phiHyper{2}{1}} }
\newcommand{\tphiz}{ {\phiHyper{2}{0}} }
\newcommand{\tPhio}{ {\PhiHyper{2}{1}} }
\newcommand{\opsio}{ {\psiHyper{1}{1}} }

%
%
\newcommand{\stirlingsubset}[2]{\genfrac{\{}{\}}{0pt}{}{#1}{#2}}
\newcommand{\stirlingcycle}[2]{\genfrac{[}{]}{0pt}{}{#1}{#2}}
\newcommand{\assocstirlingsubset}[3]{{\genfrac{\{}{\}}{0pt}{}{#1}{#2}}_{\! \ge #3}}
\newcommand{\genstirlingsubset}[4]{{\genfrac{\{}{\}}{0pt}{}{#1}{#2}}_{\! #3,#4}}
\newcommand{\irredstirlingsubset}[2]{{\genfrac{\{}{\}}{0pt}{}{#1}{#2}}^{\!\rm irr}}
\newcommand{\euler}[2]{\genfrac{\langle}{\rangle}{0pt}{}{#1}{#2}}
\newcommand{\eulergen}[3]{{\genfrac{\langle}{\rangle}{0pt}{}{#1}{#2}}_{\! #3}}
\newcommand{\eulersecond}[2]{\left\langle\!\! \euler{#1}{#2} \!\!\right\rangle}
\newcommand{\eulersecondgen}[3]{{\left\langle\!\! \euler{#1}{#2} \!\!\right\rangle}_{\! #3}}
\newcommand{\binomvert}[2]{\genfrac{\vert}{\vert}{0pt}{}{#1}{#2}}
\newcommand{\binomsquare}[2]{\genfrac{[}{]}{0pt}{}{#1}{#2}}
\newcommand{\doublebinom}[2]{\left(\!\! \binom{#1}{#2} \!\!\right)}


\newenvironment{sarray}{
             \textfont0=\scriptfont0
             \scriptfont0=\scriptscriptfont0
             \textfont1=\scriptfont1
             \scriptfont1=\scriptscriptfont1
             \textfont2=\scriptfont2
             \scriptfont2=\scriptscriptfont2
             \textfont3=\scriptfont3
             \scriptfont3=\scriptscriptfont3
           \renewcommand{\arraystretch}{0.7}
           \begin{array}{l}}{\end{array}}

\newenvironment{scarray}{
             \textfont0=\scriptfont0
             \scriptfont0=\scriptscriptfont0
             \textfont1=\scriptfont1
             \scriptfont1=\scriptscriptfont1
             \textfont2=\scriptfont2
             \scriptfont2=\scriptscriptfont2
             \textfont3=\scriptfont3
             \scriptfont3=\scriptscriptfont3
           \renewcommand{\arraystretch}{0.7}
           \begin{array}{c}}{\end{array}}


\newcommand*\circled[1]{\tikz[baseline=(char.base)]{
  \node[shape=circle,draw,inner sep=1pt] (char) {#1};}}
\newcommand{\ostar}{{\circledast}}
\newcommand{\ostarN}{{\,\circledast_{\vphantom{\dot{N}}N}\,}}
\newcommand{\ostarPsi}{{\,\circledast_{\vphantom{\dot{\Psi}}\Psi}\,}}
\newcommand{\starN}{{\,\ast_{\vphantom{\dot{N}}N}\,}}
\newcommand{\starpsi}{{\,\ast_{\vphantom{\dot{\bpsi}}\!\bpsi}\,}}
\newcommand{\starone}{{\,\ast_{\vphantom{\dot{1}}1}\,}}
\newcommand{\startwo}{{\,\ast_{\vphantom{\dot{2}}2}\,}}
\newcommand{\starinfty}{{\,\ast_{\vphantom{\dot{\infty}}\infty}\,}}
\newcommand{\starT}{{\,\ast_{\vphantom{\dot{T}}T}\,}}

\newcommand*{\Scale}[2][4]{\scalebox{#1}{$#2$}}

\newcommand*{\Scaletext}[2][4]{\scalebox{#1}{#2}} 



The hypergeometric series $\FHyper{p}{q}$ is defined by
\be
   \FHYPER{p}{q}{a_1,\ldots,a_p}{b_1,\ldots,b_q}{x}
   \;=\;
   \sum_{n=0}^\infty
   {a_1^{\overline{n}} \,\cdots\, a_p^{\overline{n}}
    \over
    b_1^{\overline{n}} \,\cdots\, b_q^{\overline{n}}
   }
   \: {x^n \over n!}
   \;,
 \label{def.pFq}
\ee
where we have used the notation $a^{\overline{n}} = a(a+1) \cdots (a+n-1)$.
In order that the series be well-defined,
we assume that $b_1,\ldots,b_q \notin -\N$
(i.e.\ no denominator parameter is a negative integer or zero);
and in order that the series not reduce to a polynomial,
we also assume that $a_1,\ldots,a_p \notin -\N$.
Then, when $p > q+1$, the series has zero radius of convergence;
when $p = q+1$, it has radius of convergence~1
and has an analytic continuation to the cut plane $\C \setminus [1,\infty)$;
and when $p \le q$, it defines an entire function of order $1/(q-p+1)$.
We are interested here in this latter case, where $\FHyper{p}{q}$ is entire.
We henceforth use the shorthand notations
$\bfa = (a_1,\ldots,a_p)$ and $\bfb = (b_1,\ldots,b_q)$,
and write $\bfa > 0$ to denote that $a_i > 0$ for all $i$
(and other similar inequalities).

A polynomial with complex coefficients is said to be
{\em negative-real-rooted}\/
if it is either identically zero or else has all its zeros in $(-\infty,0]$.
An entire function belongs to the {\em Laguerre--P\'olya class $LP^+$}\/
if it can be obtained as a limit, uniformly on compact subsets of $\C$,
of a sequence of negative-real-rooted polynomials.
Laguerre \cite{Laguerre_1882} showed in 1882 that
an entire function $f$ belongs to $LP^+$
if and only if it can be written in the form
\be
   f(x)
   \;=\;
   C x^m e^{\sigma x} \prod_{i=1}^\infty (1 + \alpha_i x)
 \label{eq.thm.LP.1}
\ee
with $C \in \C$, $m \in \N$, $\sigma,\alpha_i \ge 0$
and $\sum \alpha_i < \infty$.
See \cite[Chapter~VIII]{Levin_64} for more information
on the Laguerre--P\'olya class $LP^+$.

It is natural to investigate the conditions under which
the hypergeometric function $\FHyper{p}{q}$ ($p \le q$)
belongs to the Laguerre--P\'olya class $LP^+$.
A first result, handling the case $p=0$,
was found by Hurwitz \cite{Hurwitz_1890} in 1890,
in a paper that is unfortunately little-known;
this result was independently rediscovered by Hille \cite{Hille_29}:

\begin{theorem}[Hurwitz 1890]
   \label{thm.hurwitz}
Fix an integer $q \ge 0$.
Then for all $b_1,\ldots,b_q > 0$, the function
$\FHYPERtopzero{q}{b_1,\ldots,b_q}{\,\cdot\,}$
is an entire function of order $1/(q+1)$
that belongs to the Laguerre--P\'olya class $LP^+$.
\end{theorem}

Theorem~\ref{thm.hurwitz} is a straightforward consequence of
a (nontrivial) classical result of Laguerre
\cite[sections~16 and 17]{Laguerre_1884}
\cite[Theorem~5.6.12 and Corollary~5.6.14]{Rahman_02}.
See also \cite{Baricz_18b} for related work.

This result, along with its method of proof,
can be extended to the general case $p \le q$
under the condition that
all the parameter differences $a_i - b_i$ are nonnegative integers.
For $p=q$ this was sketched by Hille \cite{Hille_29}
and shown in detail by Ki and Kim \cite{Ki_00};
for general $p \le q$ it was obtained implicitly by
Richards \cite[pp.~477--478]{Richards_90}
and explicitly by Kalmykov and Karp \cite[Theorem~4]{Kalmykov_17}:

\begin{theorem}[Richards 1990, Kalmykov--Karp 2017]
   \label{thm.richards}
Fix integers $p \le q$.
Then for all $b_1,\ldots,b_q > 0$ and all $m_1,\ldots,m_p \in \N$,
the function
$\FHYPER{p}{q}{b_1+m_1,\ldots,b_p+m_p}{b_1,\ldots,b_q}{\,\cdot\,}$
is an entire function of order $1/(q-p+1)$
that belongs to the Laguerre--P\'olya class $LP^+$.
\end{theorem}


Finally, when $p=q$, Ki and Kim \cite[Theorem~3]{Ki_00}
proved a strong converse to Theorem~\ref{thm.richards}:

\begin{theorem}[Ki--Kim 2000]
   \label{thm.ki-kim}
Fix an integer $p \ge 1$.
If $a_1,\ldots,a_p \in \R \setminus (-\N)$ and $b_1,\ldots,b_p > 0$,
then the following are equivalent:
\begin{itemize}
   \item[(a)] $\FHyper{p}{p}(\bfa;\bfb;\,\cdot\,)$ has only a finite number
        of zeros.
   \item[(b)] $\FHyper{p}{p}(\bfa;\bfb;\,\cdot\,)$ has only real zeros.
   \item[(c)] $\FHyper{p}{p}(\bfa;\bfb;\,\cdot\,) \in LP^+$.
   \item[(d)] The $\bfa$ can be re-indexed so that $a_i = b_i + m_i$
        for $1 \le i \le p$, with all $m_i \in \N$.
\end{itemize}
\end{theorem}

It is perhaps worth remarking that in the above situation
there is an explicit formula writing $\FHyper{p}{p}(\bfa;\bfb;x)$
as $e^x$ times a polynomial of degree $|\bfm| \eqdef m_1 + \ldots + m_p$:
namely,
\be
   \FHYPER{p}{p}{b_1 + m_1 ,\,\ldots,\, b_p + m_p}
                {b_1 ,\,\ldots,\, b_p}{x}
   \;=\;
   (-1)^{|\bfm|} \,
      \Biggl( \prod\limits_{i=1}^p {1 \over b_i^{\overline{m_i}}} \Biggr)
    \: e^x \: \bfL^{(b_1 - 1,\ldots, b_p - 1)}_{m_1,\ldots,m_p}(-x)
   \;,
\ee
where $\bfL^{(\balpha)}_{\bfm}(x)$ is the (monic)
multiple Laguerre polynomial of the first kind of type~II
\cite[section~23.4.1]{Ismail_05} \cite{Sokal_multiple_laguerre}
with parameters $\balpha = (\alpha_1,\ldots,\alpha_p)$
and indices $\bfm = (m_1,\ldots,m_p)$.
When $\alpha_1,\ldots,\alpha_p > -1$
with $\alpha_i - \alpha_j \notin \Z$ for all pairs $i \neq j$,
these polynomials are multiple orthogonal \cite[Chapter~23]{Ismail_05}
with respect to the collection of measures
$x^{\alpha_i} e^{-x} \, dx$ on $(0,\infty)$ with $1 \le i \le p$;
and it follows from the general theory of multiple orthogonal polynomials
that all their zeros lie in $(0,\infty)$
\cite[Theorem~23.1.4]{Ismail_05}.
This reasoning provides an alternate proof of (d)$\implies$(c)
in Theorem~\ref{thm.ki-kim}.
(In the case $p=1$, corresponding to the ordinary Laguerre polynomials,
 this was observed long ago by Hille \cite[p.~52]{Hille_29}.)

In view of the foregoing results,
and buttressed by some numerical calculations
involving the extended Laguerre inequalities \cite{Patrick_73,Csordas_90},
Kalmykov and Karp \cite[Conjecture~3]{Kalmykov_17}
went on to conjecture that when $p < q$,
Theorem~\ref{thm.richards} could be strengthened to allow
$b_i - a_i$ to equal any positive number,
not necessarily an integer:\footnote{ 
   Kalmykov and Karp \cite[Conjecture~3]{Kalmykov_17}
   asserted only that all the zeros of $\FHyper{p}{q}(\bfa;\bfb;x)$
   are real and negative;
   but for an entire function of order $< 1$,
   this is equivalent (by Hadamard's factorization theorem)
   to being in the class $LP^+$.
   Also, Kalmykov and Karp wrote the strict inequality $a_i > b_i$;
   but if $a_i = b_i$, then the $\FHyper{p}{q}$
   trivially reduces to a $\FHyper{p-1}{q-1}$,
   so there is no harm in writing $a_i \ge b_i$.
}

\begin{conjecture}[Kalmykov--Karp 2017]
   \label{conj.kalmykov}
Suppose that $p < q$, $\bfb > 0$, and $a_i \ge b_i$ for $1 \le i \le p$.
Then the function
$\FHYPER{p}{q}{\bfa}{\bfb}{\,\cdot\,}$
is an entire function of order $1/(q-p+1)$
that belongs to the Laguerre--P\'olya class $LP^+$.
\end{conjecture}

The main purpose of the present note is to show
that Conjecture~\ref{conj.kalmykov} is false.
But what makes the situation interesting is not only
that the conjecture is false;
rather, it is as false as it can possibly be.
Namely, Theorem~\ref{thm.richards} is best possible for all $p < q$,
just as it is for $p=q$ according to Theorem~\ref{thm.ki-kim}.
We will show this by proving the following extension of
Theorem~\ref{thm.ki-kim} to $p < q$:

\begin{theorem}
   \label{thm.myconverse}
Fix integers $1 \le p \le q$.
If $a_1,\ldots,a_p \in \R \setminus (-\N)$ and $b_1,\ldots,b_p > 0$,
then the following are equivalent:
\begin{itemize}
   \item[(a)] For all $b_{p+1},\ldots,b_q > 0$,
        $\FHyper{p}{q}(\bfa;\bfb;\,\cdot\,)$ has only real zeros.
   \item[(a$\,{}'$)] For all sufficiently large $b_{p+1},\ldots,b_q > 0$,
        $\FHyper{p}{q}(\bfa;\bfb;\,\cdot\,)$ has only real zeros.
   \item[(a$\,{}''$)] For some sequence of tuples $(b_{p+1},\ldots,b_q)$
        tending to $+\infty$ in all coordinates,
        $\FHyper{p}{q}(\bfa;\bfb;\,\cdot\,)$ has only real zeros.
   \item[(b)] For all $b_{p+1},\ldots,b_q > 0$,
        $\FHyper{p}{q}(\bfa;\bfb;\,\cdot\,) \in LP^+$.
   \item[(b$\,{}'$)] For all sufficiently large $b_{p+1},\ldots,b_q > 0$,
        $\FHyper{p}{q}(\bfa;\bfb;\,\cdot\,) \in LP^+$.
   \item[(b$\,{}''$)] For some sequence of tuples $(b_{p+1},\ldots,b_q)$
        tending to $+\infty$ in all coordinates,
        $\FHyper{p}{q}(\bfa;\bfb;\,\cdot\,) \in LP^+$.
   \item[(c)] The $\bfa$ can be re-indexed so that $a_i = b_i + m_i$
        for $1 \le i \le p$, with all $m_i \in \N$.
\end{itemize}
\end{theorem}

\proof
(c)$\implies$(b) is Theorem~\ref{thm.richards},
and the implications
\begin{center}
\begin{tabular}{ccccc}
   (b) & $\Longrightarrow$ & (b${}'$) & $\Longrightarrow$ & (b${}''$) \\[2mm]
   $\big\Downarrow$  &    & $\big\Downarrow$  &    & $\big\Downarrow$ \\[2mm]
   (a) & $\Longrightarrow$ & (a${}'$) & $\Longrightarrow$ & (a${}''$)
\end{tabular}
\end{center}
\noindent
are trivial.
On the other hand, we have
\be
   \lim\limits_{b_{p+1},\ldots,b_q \to \infty}  \:
      \FHYPER{p}{q}{a_1,\ldots,a_p}{b_1,\ldots,b_q}{ b_{p+1} \cdots b_q \, x}
   \;=\;
   \FHYPER{p}{p}{a_1,\ldots,a_p}{b_1,\ldots,b_p}{x}
\ee
uniformly on compact subsets of $\C$.
So if $\FHyper{p}{q}(\bfa;\bfb;\,\cdot\,)$
has only real zeros (or is in $LP^+$)\footnote{
   In fact, for an entire function of order $< 1$
   (as is the case when $p < q$),
   having only real nonpositive zeros is {\em equivalent}\/ to being in $LP^+$.
   And a hypergeometric function with $p \le q$ and parameters $\bfa,\bfb > 0$
   obviously cannot have positive real zeros.

   On the other hand, when the entire function is of order~1
   (as it is when $p=q$), this equivalence does not hold,
   so the equivalence of (b) and (c) in Theorem~\ref{thm.ki-kim}
   is a nontrivial fact.
}
for some sequence of tuples $(b_{p+1},\ldots,b_q)$
tending to $+\infty$ in all coordinates,
then the same holds for $\FHyper{p}{p}(\bfa;b_1,\ldots,b_p;\,\cdot\,)$.
In this situation Theorem~\ref{thm.ki-kim} implies that
the $\bfa$ can be re-indexed so that $a_i = b_i + m_i$;
so (a${}''$) or (b${}''$) implies (c).
\qed

In other words, Theorem~\ref{thm.myconverse}
is an almost trivial corollary of Theorem~\ref{thm.ki-kim}.
But since the proof \cite{Ki_00} of Theorem~\ref{thm.ki-kim}
is far from trivial, it is also of some value to exhibit
explicit counterexamples to Conjecture~\ref{conj.kalmykov};
we do this in the Appendix.

Of course, Theorem~\ref{thm.myconverse} does not exclude that,
for noninteger values of $b_i - a_i$ ($1 \le i \le p$),
the function $\FHyper{p}{q}(\bfa;\bfb;\,\cdot\,)$
can belong to $LP^+$ for some bounded range of values of $b_{p+1},\ldots,b_q$.
For instance, Driver {\em et al.}\/ \cite[proof of Theorem~8]{Driver_07}
have observed that
\be
   \FHYPER{1}{2}{a-\smhalf}{a,\, 2a-1}{x}
   \;=\;
   \left[ \FHYPERtopzero{1}{a}{x/4} \right]^{2}
\ee
and hence (by Theorem~\ref{thm.hurwitz}) that
$\displaystyle \FHYPER{1}{2}{a-\smhalf}{a,\, 2a-1}{ \,\cdot\, } \in LP^+$
for all $a > 0$.
More generally, they observe that \cite[eq.~(2.03)]{Bailey_28}
\be
   \FHYPER{2}{3}{(a+b)/2,\,(a+b-1)/2}{a,\,b,\,a+b-1}{x}
   \;=\;
   \FHYPERtopzero{1}{a}{x/4} \; \FHYPERtopzero{1}{b}{x/4}
\ee
and hence that
$\displaystyle \FHYPER{2}{3}{(a+b)/2,\,(a+b-1)/2}{a,\,b,\,a+b-1}{ \,\cdot\, }
 \in LP^+$ for all $a,b > 0$.
Likewise, P\'olya \cite[p.~379]{Polya_18} and Hille \cite[p.~53]{Hille_29}
have shown that
$\displaystyle \FHYPER{1}{2}{1}{a,\, a+\smhalf}{ \,\cdot\, } \in LP^+$
for $0 < a \le 3/2$
(but has no real zeros when $a > 3/2$);
see also \cite[Theorem~3.6]{Craven_06} for the special case $a = 1/4$,
and \cite{Cho_18,Cho_20} for further cases of the absence of real zeros.
Finally, Craven and Csordas \cite[Proposition~3.11 and Theorem~3.13]{Craven_06}
have shown that
$\displaystyle
 \FHYPER{1}{2}{\smhalf}{\smfrac{1}{3},\, \smfrac{2}{3}}{ \,\cdot\, } \in LP^+$.
(See also \cite[pp.~219--220]{Baricz_18a} for a $\FHyper{3}{4}$
 that belongs to $LP^+$.)
It~would be interesting to determine the exact set of parameters
$(a_1,b_1,b_2) \in \R^3$ for which the function
$\FHyper{1}{2}(a_1; b_1,b_2; \,\cdot\,)$ belongs to $LP^+$;
but this is probably a hopeless task,
as the boundary of this set is probably not given by any simple formula.
If that is the case, then the best we can do
is to find decent inclusion or exclusion regions,
or monotonicities guaranteeing that
$\FHyper{1}{2}(a_1; b_1,b_2; \,\cdot\,) \in LP^+$
implies
$\FHyper{1}{2}(a'_1; b'_1,b'_2; \,\cdot\,) \in LP^+$
under suitable conditions relating
$(a_1,b_1,b_2)$ to $(a'_1,b'_1,b'_2)$.
Indeed, I~am unaware of any 2-dimensional or 3-dimensional sets
of $(a_1,b_1,b_2) \in \R^3$ for which $\FHyper{1}{2}$
has been proven to belong to $LP^+$.
These are interesting tasks for future research,
which are complementary to work \cite{Cho_18,Cho_20}
on the absence of real zeros.

\appendix
\section{Explicit counterexamples to Conjecture~\ref{conj.kalmykov}}

In order to exhibit counterexamples to Conjecture~\ref{conj.kalmykov},
we proceed as follows.
To test whether a function $f$ belongs to $LP^+$,
we will use the following criterion\footnote{
   Results similar to Proposition~\ref{prop.criterion} go back at least
   to the work of Grommer \cite[especially pp.~157--158]{Grommer_14} in 1914.
   See also the article of Krein \cite{Krein_62} for a very useful survey.
}:

\begin{proposition}[Logarithmic-derivative criterion for $f \in LP^+$]
   \label{prop.criterion}
For an entire function $f$ with $f(0) \neq 0$, the following are equivalent:
\begin{itemize}
   \item[(a)]  $f \in LP^+$.
   \item[(b)]  The sequence
$\displaystyle \Bigl( (-1)^n \, [x^n] {f'(x) \over f(x)} \Bigr)_{n \ge 0}$
is a Stieltjes moment sequence.
   \item[(c)]  The sequence
$\displaystyle \Bigl( (-1)^n \, [x^n] {f'(x) \over f(x)} \Bigr)_{n \ge 0}$
is a Stieltjes moment sequence with a unique representing measure
of the special form
\be
   \mu  \;=\;  \sigma \delta_0 \:+\: \sum_i m_i \alpha_i \, \delta_{\alpha_i}
 \label{eq.prop.criterion}
\ee
for a sequence $\alpha_1 > \alpha_2 > \ldots > 0$
satisfying $\sum\limits_i \alpha_i < \infty$,
integers $m_i \ge 0$, and $\sigma \ge 0$.
\end{itemize}
\end{proposition}

\noindent
Here we write $[x^n] \, g(x)$ to denote the coefficient of $x^n$
in the Taylor expansion of $g(x)$;
and we recall that a sequence is called a {\em Stieltjes moment sequence}\/
if it is the moment sequence for some positive measure on $[0,\infty)$.
The proof of Proposition~\ref{prop.criterion} is not difficult:
logarithmic differentiation of \reff{eq.thm.LP.1}
leads to
\be
   {f'(x) \over f(x)}
   \;=\;
   \sigma \:+\: \sum_{i=1}^\infty {\alpha_i \over 1 + \alpha_i x}
   \;,
\ee
which shows that (a)$\iff$(c);
(c)$\implies$(b) is trivial;
and a short argument shows that (b)$\implies$(c) holds
because of the assumption that $f$ is entire.
In what follows we will use only the elementary fact (a)$\implies$(b).

As a test for when a sequence of real numbers is a Stieltjes moment sequence,
we will use the following criterion,
due to Stieltjes \cite{Stieltjes_1894} in 1894:

\begin{proposition}[Continued-fraction criterion for Stieltjes moment property]
   \label{prop.stieltjes}
\hfill\break
A sequence $\ba = (a_n)_{n \ge 0}$ of real numbers
is a Stieltjes moment sequence if and only if its ordinary generating function
has a continued-fraction expansion
\be
   \sum_{n=0}^{\infty} a_n t^n
   \;=\;
   \cfrac{\alpha_0}{1 - \cfrac{\alpha_1 t}{1 - \cfrac{\alpha_2 t}{1 - \cdots}}}
   \label{eq.thm.hankelreal.infty.Stype}
\ee
in the sense of formal power series,
with {\em nonnegative}\/ coefficients $\alpha_1,\alpha_2,\ldots\;$.
\end{proposition}

We refer to \cite{Shohat_43,Akhiezer_65,Simon_98,Schmudgen_17}
for further information on the moment problem.

\bigskip

We will use the combination of
Propositions~\ref{prop.criterion} and \ref{prop.stieltjes}
to show that certain entire functions $f$ do not belong to $LP^+$.
Please observe that in order to show, using Proposition~\ref{prop.stieltjes},
that a given sequence $\ba$ is not Stieltjes,
we will need to show not only that some $\alpha_n < 0$
but also that $\alpha_1,\ldots,\alpha_{n-1} \neq 0$
(so that the continued fraction does not terminate prematurely).
But in our applications each $\alpha_i$ will vanish (if at all)
on a subvariety of codimension~1 in the parameter space,
so this vanishing will be harmless if we can show that $\alpha_n < 0$
on a nonempty open set.

We will consider here the simplest case of Conjecture~\ref{conj.kalmykov}
that is not contained in Theorem~\ref{thm.hurwitz},
namely $p=1$ and $q=2$.
The logarithmic derivative of $f(x) = \FHYPER{1}{2}{a_1}{b_1,b_2}{x}$ is
\begin{eqnarray}
   & &
   {f'(x) \over f(x)}
   \;=\;
   {a_1 \over b_1 b_2}
     \:-\: {a_1 [a_1 (1 + b_1 + b_2) - b_1 b_2]  
            \over
            b_1^2 (b_1 + 1) b_2^2 (b_2 + 1)
           }  \: t
      \nonumber \\[1mm]
   & & \hspace{2.65cm}
     \:+\: {\hbox{a rather complicated expression}
            \over
            b_1^3 (b_1 + 1) (b_1 + 2) b_2^3 (b_2 + 1) (b-2 + 2)
           } \: t^2
     \:-\: \ldots
   \qquad
\end{eqnarray}
We divide this series by the leading coefficient $a_1/(b_1 b_2)$,
introduce the signs $(-1)^n$,
and then use the Euler--Viscovatov algorithm
\cite{Euler_1760,Viscovatov_1806} \cite[pp.~27--31]{Khovanskii_63}
\cite{Sokal_alg_contfrac}
to compute the continued-fraction coefficients $\alpha_1,\alpha_2,\ldots\;$.
They are rational functions
\be
   \alpha_i(a_1,b_1,b_2)  \;=\;  {N_i(a_1,b_1,b_2) \over D_i(a_1,b_1,b_2)}
   \;,
\ee
and the first two are
\begin{subeqnarray}
   \alpha_1  & = &  {a_1 (1 + b_1 + b_2) - b_1 b_2
                     \over
                     b_1 (b_1 + 1) b_2 (b_2 + 1)
                    }
            \\[2mm]
   \alpha_2  & = &
      {(1 + a_1) \bigl[ a_1 (2 + 3b_1 + 3b_2 + b_1^2 + b_1 b_2 + b_2^2)
                     \,-\, b_1 b_2 (3 + b_1 + b_2) \bigr]
       \over
       (b_1 + 1) (b_1 + 2) (b_2 + 1) (b_2 + 2)
          \bigl[ a_1 (1 + b_1 + b_2) - b_1 b_2 \bigr]
      }
   \qquad
\end{subeqnarray}
The subsequent $\alpha_i$ become increasingly complicated.

It is convenient to write $a_1 = b_1 + \gamma$.
Then Conjecture~\ref{conj.kalmykov} states that, for all $i$,
$\alpha_i \ge 0$ whenever $b_1, b_2 > 0$ and $\gamma \ge 0$.
We will now show that this is false, by looking at $\alpha_3$.
The denominator polynomial
\begin{eqnarray}
   & &
   D_3
   \;=\;  
   (b_1 + 2) (b_1 + 3) (b_2 + 2) (b_2 + 3)
   \bigl[ b_1 (b_1 + 1) \,+\, (1 + b_1 + b_2) \gamma \bigr]
      \:\times\:
      \nonumber \\
   & & \qquad\qquad
   \bigl[ b_1 (b_1 + 1) (b_1 + 2) \,+\,
          [(b_1 + 1)(b_1 + 2) + (3 + b_1 + b_2) b_2] \gamma \bigr]
\end{eqnarray}
has nonnegative coefficients and is therefore strictly positive whenever
$b_1, b_2 > 0$ and $\gamma \ge 0$.
The numerator polynomial $N_3$ is a polynomial of degree~4 in $b_2$,
and its highest-order term is
\be
   [b_2^4] \, N_3
   \;=\;
   \gamma^2 (\gamma - 1)
   \;,
\ee
which is negative whenever $0 < \gamma < 1$
(and also when $\gamma < 0$).
It follows that, whenever $b_1 > 0$ and $0 < \gamma < 1$ (or $\gamma < 0$),
we have $N_3 < 0$ and hence $\alpha_3 < 0$
for all sufficiently large positive $b_2$.
For instance, if $b_1 = 1$ and $\gamma = 1/2$,
then $\alpha_3 < 0$ for all $b_2 > 52.4865\ldots\;$.
By Propositions~\ref{prop.criterion} and \ref{prop.stieltjes},
it follows that $\FHyper{1}{2}(b_1 + \gamma; b_1, b_2; \,\cdot\,) \notin LP^+$
in these cases.

Though the $LP^+$ property fails for $0 < \gamma < 1$,
one might hope that it is restored when $\gamma \ge 1$.
But this too is false, as we can see by looking at $\alpha_5$.
One can verify, using {\sc Mathematica}'s function {\tt Reduce},
that the denominator polynomial $D_5$ is strictly positive
whenever $b_1, b_2 \ge 0$ and $\gamma \ge 1$.
(Unlike what occurred for $D_3$, it is not {\em coefficientwise}\/ nonnegative
 in $b_1$, $b_2$ and $\gamma-1$.)
On the other hand, the numerator polynomial $N_5$
is a polynomial of degree~11 in $b_2$,
and its highest-order term is
\be
   [b_2^{11}] \, N_5
   \;=\;
   \gamma^4 (\gamma - 1)^2 (\gamma - 2)
   \;,
\ee
which is negative whenever $1 < \gamma < 2$
(and also when $\gamma < 0$ or $0 < \gamma < 1$).
It follows that, whenever $b_1 > 0$ and $1 < \gamma < 2$,
we have $N_5 < 0$ and hence $\alpha_5 < 0$
for all sufficiently large positive $b_2$.
For instance, if $b_1 = 1$ and $\gamma = 3/2$,
then $\alpha_5 < 0$ for all $b_2 > 574.8859\ldots\;$.

This pattern appears to continue.
Using {\sc Mathematica}'s function {\tt FindInstance},
I found (and the reader can easily verify) that
$\alpha_7 < 0$ when $b_1 = 1$, $\gamma = 5/2$ and $b_2 = 72053$,
and that
$\alpha_9 < 0$ when $b_1 = 1$, $\gamma = 7/2$ and $b_2 = 750232$.
More precisely, it appears that:

\begin{conjecture}[Sufficient conditions for $\FHyper{1}{2} \notin LP^+$]
   \label{conj.1F2}
\hfill\break
\vspace*{-6mm}
\begin{itemize}
   \item[(a)]  The denominator polynomial $D_n$ is strictly positive
       whenever $b_1, b_2 > 0$ and $\gamma \ge \lfloor (n-2)/2 \rfloor$.
   \item[(b)]  The numerator polynomial $N_n$ is of degree~$\binom{n}{2} + 1$
       in $b_2$, and its highest-order coefficient is
\be
   [b_2^{\binom{n}{2} + 1}] \, N_n
   \;=\;
   \begin{cases}
       (\gamma-k) \prod\limits_{i=0}^{k-1} (\gamma-i)^{2k-2i}
                      & \textrm{if $n=2k+1$}  \\[2mm]
       (b_1 + \gamma + k+1)
                  \prod\limits_{i=0}^{k} (\gamma-i)^{2k-2i+1}
                      & \textrm{if $n=2k+2$}
   \end{cases}
\ee
\end{itemize}
and hence
\begin{itemize}
   \item[(c)] For any $b_1 > 0$ and any $\gamma \in (k-1,k)$,
      we have $\alpha_{2k+1} < 0$ for all sufficiently large $b_2$.
   \item[(d)] For any $b_1 > 0$ and any noninteger $\gamma > 0$,
      we have $\FHyper{1}{2}(b_1 + \gamma; b_1,b_2; \,\cdot\,) \notin LP^+$
      for all sufficiently large $b_2$.
\end{itemize}
\end{conjecture}

I have checked part~(a) for $1 \le n \le 5$,
and part~(b) for $1 \le n \le 10$.

\section*{Acknowledgments}

I wish to thank Kathy Driver for drawing my attention to the work of
Hurwitz \cite{Hurwitz_1890} and Hille \cite{Hille_29},
Dmitrii Karp for drawing my attention to \cite{Baricz_18a},
and Alex Dyachenko for helpful comments.

This research was supported in part by
the U.K.~Engineering and Physical Sciences Research Council grant EP/N025636/1.

\end{document}